\documentclass[10pt]{article}

\usepackage{latexsym,amsfonts}
\usepackage{amssymb,amsthm,upref,amscd}
\usepackage[T1]{fontenc}
\usepackage{times}

\newtheorem{Theorem}{Theorem}[section]
\newtheorem{definition}[Theorem]{Definition}
\newtheorem{lemma}[Theorem]{Lemma}

\newtheorem{corollary}[Theorem]{Corollary}
\newtheorem{remark}[Theorem]{Remark}
\newtheorem{Open Problem}[Theorem]{Open Problem}

\makeatletter \@addtoreset{equation}{section} \makeatother
\begin{document}

\title{\bf Nontrivial solutions for periodic  Schr\"odinger
equations with sign-changing nonlinearities
 }

\author{Shaowei Chen \thanks{ E-mail address:
swchen6@163.com (Shaowei Chen)}    \quad \quad Conglei Wang
   \quad \quad Liqin Xiao\\ \\
\small \small\it School of Mathematical Sciences, Huaqiao
University, \\
\small Quanzhou  362021,  China\\ }

\date{}
\maketitle
\begin{minipage}{13cm}
{\small {\bf Abstract:} Using a new infinite-dimensional linking
theorem, we obtained  nontrivial solutions for  strongly
indefinite periodic
Schr\"odinger equations with sign-changing nonlinearities. 
}
\end{minipage}

\section{Introduction and statement of results}\label{diyizhang}

In this paper,  we consider the following semilinear
Schr\"{o}dinger equation:
\begin{equation}\label{iiiuyhgt}
- \Delta u+V(x)u+f(x,u)=0,\qquad u\in H^{1}(\mathbb{R}^{N}),
\end{equation}
where $N\geq1$. The potential $V$ is a continuous function and
1-periodic in $x_j$ for $j = 1,\cdots,N$. In this case, the
spectrum of the operator $-\Delta+V$ is a purely continuous
spectrum that is bounded below and consists of closed disjoint
intervals (\cite[Theorem XIII.100]{Reed2}). Thus, the complement
$\mathbb{R} \setminus \sigma(L)$ consists of  open intervals
called spectral gaps. More precisely, for $V$, we assume
\begin{description}
\item{$(\bf{v}).$}
 $V\in C( \mathbb{R}^{N}) $ is 1-periodic in $x_j$ for $j = 1,\cdots,N$,
  0 is in a spectral gap $(-\mu_{-1}, \mu_1)$ of  $-\Delta+V$
  and $-\mu_{-1}$ and $\mu_1$ lie in the essential spectrum of $-\Delta+V.$
Denote $$\mu_0 := \min\{\mu_{-1}, \mu_1\}.$$
\end{description}
And for $f,$ we assume
\begin{description}
\item{$(\bf{f_1}).$} $f\in C( \mathbb{R}^{N}\times\mathbb{R}) $ is
1-periodic in $x_j$ for $j = 1,\cdots,N$.  And there exist
constants $C>0$ and $2< p< 2^*$ such that
$$|f(x,t)|\leq C(1+|t|^{p-1}),\ \forall (x,t)\in\mathbb{R}^N\times\mathbb{R}$$
where $2^*:= \left\{
\begin{array}{l}
\frac{2N}{N-2}, \ N\geq 3\\
\infty,  \quad N=1,2.\\
\end{array} \right.$

\item{$(\bf{f_2}).$} The limits $\lim_{t\rightarrow 0}f( x,t)/t =
0 $ and $\lim_{|t|\rightarrow \infty}f( x,t)/t = +\infty $ hold
uniformly for $x\in\mathbb{R}^{N}$.
\end{description}

Let  $F(x,t)=\int^{t}_{0}f(x,s)ds$, $V_-(x)=\max\{-V(x), 0\}$,
 $x\in\mathbb{R}^N,$ and
\begin{eqnarray}\label{mmmmauuhdff7}
\widetilde{F}(x,t):=\frac{1}{2}tf(x,t)-F(x,t).
\end{eqnarray}
\begin{description}
\item{$(\bf{f_3}).$} There exist $\rho>0$ and $M>\rho$ such that
\begin{eqnarray}\label{nc99ciduff}
\inf_{|t|\geq M/2,\
x\in\mathbb{R}^N}\Big(\frac{f(x,t)}{t}-V_-(x)\Big)>0,
\end{eqnarray}
\begin{eqnarray}\label{nc99cidudyff}
\inf_{|t|\geq\rho, x\in\mathbb{R}^N}F(x,t)\geq 0,
\end{eqnarray}

\begin{eqnarray}\label{8r7bcvgftfudyff}
\inf_{\rho\leq|t|\leq M, x\in\mathbb{R}^N}\widetilde{F}(x,t)> 0,
\end{eqnarray}

and
 \begin{eqnarray}\label{l99ds7yydt}
 \sup_{|t|\leq\rho, x\in\mathbb{R}^N}\Big|\frac{f(x,t)}{t}\Big|
 +\Big(\sup_{\rho\leq|t|\leq M, x\in\mathbb{R}^N}\frac{f^2(x,t)}{|\widetilde{F}(x,t)|}\Big)^{\frac{1}{2}}\cdot
 \Big(\sup_{|t|\leq\rho,
 x\in\mathbb{R}^N}\frac{|\widetilde{F}(x,t)|}{t^2}\Big)^{\frac{1}{2}}
 <\mu_0.
 \end{eqnarray}
  \end{description}

\medskip

A solution $u$ of (\ref{iiiuyhgt}) is called nontrivial if
$u\not\equiv0.$ Our main result is the following theorem:

\begin{Theorem}\label{th1}
Suppose that $\bf ( v) $ and $\bf (f_1)-(f_3)$ are satisfied. Then
the problem (\ref{iiiuyhgt}) has a nontrivial solution.
\end{Theorem}

As an application of  this theorem, we have the following
corollary:

\begin{corollary}\label{nc88fyyyf2w}
Suppose that $2<p<q<2^*,$ and $\bf ( v) $ is satisfied. Then there
exists $\lambda_0>0$ such that, for  $0<\lambda<\lambda_0,$
$$-\Delta u+Vu=\lambda|u|^{p-2}u-|u|^{q-2}u$$ has a nontrivial solution.
\end{corollary}
\begin{remark}\label{kkkkkk8uyyta}
This corollary shows that under the assumptions $\bf (f_1)-(f_3)$,
the nonlinearity $f$  may be sign-changing.
\end{remark}

Under the assumption $\bf (v)$, the quadratic form
$\int_{\mathbb{R}^N}(|\nabla u|^2+V(x)u^2)dx$ has
infinite-dimensional negative and positive spaces. This case is
called strongly indefinite. Semilinear periodic Schr\"odinger
equations with strongly indefinite linear part have attracted much
attention in recent
 years due to its numerous  applications in mathematical physics.
In \cite{alama}, the authors used a dual variational method to
obtain a nontrivial solution of (\ref{iiiuyhgt}) with
$f(x,t)=-W(x)|t|^{p-2}t$, where $W$ is a perturbed  periodic
function and $2<p<2^*.$  In \cite{TW}, Troestler  and Willem used
critical point theory  to obtain  a nontrivial solution of
(\ref{iiiuyhgt}) by assuming that $f\in C^1(\mathbb{R}^N\times
\mathbb{R})$, $|\partial_tf(x,t)|\leq C(|t|^{p-2}+|t|^{q-2})$ for
some $2<p<q<2^*$ and $g=-f$ satisfies the so-called
Ambrosetti-Rabinowitz condition $0<\gamma G(x,t)\leq tg(x,t)$,
$\forall (x,t)\in \mathbb{R}^N\times (\mathbb{R}\setminus\{0\})$,
where $G(x,t)=\int^t_0g(x,s)ds.$  Then, Kryszewski and Szulkin
\cite{KS} proved a new infinite-dimensional linking theorem. Using
it, they generalized Troestler and Willem's result by assuming
that $f\in C(\mathbb{R}^N\times \mathbb{R})$ and $g=-f$ satisfies
the
 Ambrosetti-Rabinowitz condition.
Similar results were also obtained by Pankov and Pfl\"uger in
\cite{pankov, PMilan} by  an approximation method and a variant
Nehari method.  Equation (\ref{iiiuyhgt}) with asymptotically
linear nonlinearities and other super-linear nonlinearities has
also been studied by many authors. One can see \cite{chenzhang,
Dingbook, yanheng, LiSzulkin} for the asymptotically linear case
and \cite{Ackermann, BC, chen, yanheng, KS1, liu, SW, Y} for the
super-linear case. Moreover, equation (\ref{iiiuyhgt}) with $0$
belonging to the spectrum of $-\Delta+V$ was studied in
\cite{Dingbook, WZ, YCD}. Finally, we should mention that the
methods of studying the strongly indefinite periodic Schr\"odinger
can shed light  on  other strongly indefinite problems, such as,
the Hamiltonian systems or the discrete nonlinear Schr\"odinger
equations. One can consult \cite{ Dingbook, szulkin} or
\cite{pankovJMAA, zhangguop}.

 All the existence results for equation (\ref{iiiuyhgt}) we
 mentioned above
 are obtained under the assumption that $f$ does not change sign in  $\mathbb{R}^N\times \mathbb{R}$, i.e.,
 $f\geq 0$ in $\mathbb{R}^N\times \mathbb{R}$ or $f\leq 0$ in $\mathbb{R}^N\times
 \mathbb{R}$.
 However, under our  assumptions $\bf (f_1)-(f_3)$,
 $f(x,t)$ can be negative in
 $\{(x,t)\ |\ |t|<\rho\}$. Together with $(\ref{nc99cidudyff})$, this implies that
  $f(x,\cdot)$ may change sign in $ \mathbb{R}$. As we know, this
  situation has never been studied before. And it is the novelty
  of our main results Theorem \ref{th1} and Corollary \ref{nc88fyyyf2w}.

  The difficulties  of equation (\ref{iiiuyhgt}) with sign-changing
  nonlinearity come from two aspects. The first  is that the
  classical infinite-dimensional linking theorem (see
  \cite[Theorem 6.10]{willem} or \cite{KS}) cannot   be used to deal (\ref{iiiuyhgt}) in this
  case. To use this
  linking theorem, the functional
  corresponding to (\ref{iiiuyhgt}) must satisfy some upper
  semi-continuous assumption. However, when $f(x,t)$ is sign-changing,
  the functional
  corresponding to (\ref{iiiuyhgt}) does not satisfy this
  assumption. The second is that the sign-changing nonlinearity
  brings more difficulty in the proof of boundedness of
  Palais-Smale sequence.

  In this paper,  a variant  infinite-dimensional linking
  theorem (see Theorem \ref{b99d443}) is given. This theorem replaced the upper semi-continuous
  assumption in the
  classical infinite-dimensional linking theorem (\cite[Theorem 6.10]{willem}) with other assumptions.
  For the reader's convenience, we presented this theorem and its proof in the appendix.
Using this theorem, a $(\overline{C})_c$ sequence (see Definition
\ref{bcvujnbytg}) of equation (\ref{iiiuyhgt}) is obtained. Under
$\bf (f_1)-(f_3)$, we can prove that this sequence is bounded in
$H^1(\mathbb{R}^N)$ (see Lemma \ref{x4esdftrg}). Then, we obtained
a nontrivial solution of (\ref{iiiuyhgt}) from the
$(\overline{C})_c$ sequence  through the concentration-compactness
principle.

 \medskip

 \noindent{\bf Notation.} $B_r(a)$ denotes the  open ball of radius $r$ and center $a$.
For a Banach space $E,$ we denote the dual space of $E$ by $E'$,
and denote   strong and   weak convergence in $E$  by
$\rightarrow$ and $\rightharpoonup$, respectively. For $\varphi\in
C^1(E;\mathbb{R}),$ we denote the Fr\'echet derivative of
$\varphi$ at $u$ by $\varphi'(u)$.  The Gateaux derivative of
$\varphi$ is denoted by $\langle \varphi'(u), v\rangle,$ $\forall
u,v\in E.$ $L^p(\mathbb{R}^N)$ denotes the standard $L^p$ space
$(1\leq p\leq\infty)$, and $H^1(\mathbb{R}^N)$ denotes the
standard Sobolev space with norm
$||u||_{H^1}=(\int_{\mathbb{R}^N}(|\nabla u|^2+u^2)dx)^{1/2}.$  We
use $O(h)$, $o(h)$ to mean $|O(h)|\leq C|h|$ and
$o(h)/|h|\rightarrow 0$.

\section{Variational  and linking structure for equation (\ref{iiiuyhgt}) }\label{mc88vuyfyf}

Under the  assumptions $\bf (v)$,  and $\bf (f_1)$, the functional
\begin{equation}\label{mnncb66dtd}
J(u)=\frac{1}{2}\int_{\mathbb{R}^{N}}\left\vert \nabla
u\right\vert ^{2}dx
+\frac{1}{2}\int_{\mathbb{R}^{N}}V(x)u^{2}dx+\int
_{\mathbb{R}^{N}}F(x,u)dx
\end{equation}
is of class $C^1$ on $X := H^1(\mathbb{R}^N )$. The derivative of
$J$ is
\begin{eqnarray}\label{nv88ufyyfh1}
\langle J'(u), v\rangle=\int_{\mathbb{R}^N}(\nabla u\nabla
v+V(x)uv)dx+\int_{\mathbb{R}^N}f(x,u)vdx,\ u,v\in
H^1(\mathbb{R}^N)
\end{eqnarray} and the
critical points of $J$ are weak solutions of (\ref{iiiuyhgt}).

There is a standard  variational setting for the quadratic form
$\int_{\mathbb{R}^N}(|\nabla u|^2+V(x)u^2)dx$. For the reader's
convenience, we state it here.  One can consult \cite{yanheng} or
\cite{Dingbook} for more details.

Assume that $\bf (v)$ holds and let $S=-\Delta+V$ be the
self-adjoint operator acting on $L^2(\mathbb{R}^N)$ with domain
$D(S)=H^2(\mathbb{R}^N)$. By virtue of $\bf (v)$,  we have the
orthogonal decomposition $$L^2=L^2(\mathbb{R}^N)=L_1+L_2$$ such
that $S$ is   positive (resp.  negative) in $L_1$(resp.in $L_2$).
Let $X=D(|S|^{1/2})$ be equipped with the inner product
$$(u,v)=(|S|^{1/2}u, |S|^{1/2}v)_{L^2}$$and norm $||u||=|||S|^{1/2}u||_{L^2}$,
 where $ (\cdot,\cdot)_{L^2}$ denotes the
inner product of $L^2$. From $\bf (v),$ $$X=H^1(\mathbb{R}^N)$$
with equivalent norms. Therefore, $X$  continuously embeds in $L^q
(\mathbb{R}^N)$ for all $2\leq q\leq 2^*$. In addition, we have
the decomposition $$X=Y\oplus Z,$$ where $Y=X\cap L_1$ , $Z=X\cap
L_2$ and $Y,$ $Z$ are orthogonal with respect to both
$(\cdot,\cdot)_{L^2}$ and $(\cdot,\cdot)$. Let  $P:X\rightarrow Y$
and $Q:X\rightarrow Z$ be orthogonal projections.  Therefore, for
every $u\in X$ , there is a unique decomposition
$$ u=Pu+Qu$$
 with $(Pu,Qu)=0$
and \begin{eqnarray}\label{nbvbviifjfjj}\int_{\mathbb{R}^N}|\nabla
u|^2dx+\int_{\mathbb{R}^N}V(x)u^2dx=||Pu||^2-||Qu||^2,\ u\in X.
\end{eqnarray}
 Moreover,
\begin{eqnarray}\label{bcvttdref}
&&\mu_{-1}||Pu||^2_{L^2}\leq||Pu||^2,\quad \forall u\in X,
\end{eqnarray} and\begin{eqnarray}\label{bcvttdref2}
&&\mu_{1}||Qu||^2_{L^2}\leq||Qu||^2,\quad \forall u\in X.
\end{eqnarray}
Therefore, \begin{eqnarray}\label{mmnc999ciu}
\mu_0||u||^2_{L^2}\leq ||u||^2,\ \forall u\in X.
\end{eqnarray}
We denote
$$X^+:=Z,\ X^-:=Y,\ u^+:=Qu,\ u^-:=Pu.$$
Let $$\Phi(u)=-J(u),\ u\in X.$$  Then by (\ref{mnncb66dtd}) and
(\ref{nbvbviifjfjj}), $\Phi$ can be written as
\begin{eqnarray}\label{vcrdfd}
 \Phi(u)=\frac{1}{2}(||u^+||^2-||u^-||^2)-\int_{\mathbb{R}^N}F(x,u)dx,\ u\in X.
\end{eqnarray}
And the derivative of $\Phi$ is given by
\begin{eqnarray}\label{nc88c766f5}
\langle\Phi'(u),
v\rangle=(u^+,v)-(u^-,v)-\int_{\mathbb{R}^N}f(x,u)v dx,\ \forall
u,v\in X.
\end{eqnarray}

Let $\{e^-_k\}$ be the total orthonormal sequence in $X^-$.   We
define
\begin{eqnarray}\label{o88dytgdg}
|||u|||=\max\Big\{||Qu||,\sum^{\infty}_{k=1}\frac{1}{2^{k+1}}|(Pu,e^-_k)|\Big\}
\end{eqnarray}
on $X.$

The following Lemma shows that $\Phi$ satisfies the linking
condition (see (\ref{nx99s8s7yy})) of Theorem \ref{b99d443} in the
appendix.

\begin{lemma}\label{bc77ftrg} Suppose $\bf (v)$,  and $\bf ( f_{1})-\bf ( f_{3}) $  are satisfied.
Then there exist $\delta>0,$  $R > r > 0$ and  $u_0 \in X^+$ with
$||u_0||
 = 1$ such that
 \begin{eqnarray}\label{m99a66tsrd}
 \inf_N\Phi>\max\Big\{ \sup_{\partial M}\Phi,
 \sup_{|||u|||\leq\delta}\Phi\Big\},
 \end{eqnarray}
 where $N,$ $M$ and $\partial M$ are defined in (\ref{ncg7777qwwq})
 and (\ref{m9155sfddd}) in the appendix.
\end{lemma}
\noindent{\bf Proof.} We divide the proof into several steps.

{\bf Step 1.} We shall prove that there exists $r>0$ such that
$\inf_N\Phi>0$.

 From $\bf (f_1)$ and $\bf (f_2)$, we deduce that for any
 $\epsilon>0,$ there exists $C_\epsilon>0$ such that
 \begin{eqnarray}\label{mn88fyfgg}
 |F(x,t)|\leq\epsilon t^2+C_\epsilon|t|^p,\ \forall
 t\in\mathbb{R}.\nonumber
 \end{eqnarray}
Then by the Sobolev inequality $||u||_{L^p(\mathbb{R}^N)}\leq C
||u||$, $\forall u\in X$ and the definition of $\Phi$ (see
(\ref{vcrdfd})), there exists $C'>0$ such that, for any $u\in
X^+,$
\begin{eqnarray}\label{nn7cxrxex}
\Phi(u)\geq \frac{1}{2}||u||^2-C'\epsilon
||u||^2-C'C_\epsilon||u||^p=(\frac{1}{2}-C'\varepsilon)||u||^2-C'C_\epsilon||u||^p.
\end{eqnarray}
Choose  $\epsilon=1/4C'$ in  (\ref{nn7cxrxex}) and let
$r=(8C'C_\epsilon)^{-\frac{1}{p-2}}$. We get that, for $N=\{u\in
X^+\ |\ ||u||=r\}$,
\begin{eqnarray}
\inf_N \Phi\geq r^2/8>0.\nonumber
\end{eqnarray}

{\bf Step 2.} We shall prove that $\Phi(u)\rightarrow -\infty$ as
$||u||\rightarrow\infty$ and $u\in X^-\oplus\mathbb{R}^+u_0$. As a
consequence, there exists $R>r$ such that $\sup_{\partial
M}\Phi\leq0<\inf_N \Phi.$

Arguing indirectly, assume that for some sequence
 $u_n\in X^-\oplus \mathbb{R}^+u_0$ with
$||u_n||\rightarrow+\infty$, there is $ \mathcal{L}
> 0$ such that $\Phi(u_n)\geq-\mathcal{L}$ for all $n$. Then, setting
$w_n=u_n/||u_n||$,   we have $||w_n||=1$, and, up to a
subsequence, $w_n\rightharpoonup w$, $w^-_n\rightharpoonup w^-\in
X^-$ and $w^+_n\rightarrow w^+\in X^+.$
 Dividing both sides of
$\Phi(u_n)\geq-\mathcal{L}$ by $||u_n||^2$, we get that
\begin{eqnarray}\label{mmvn77dtdg}
-\frac{\mathcal{L}}{||u_n||^2}\leq\frac{\Phi(u_n)}{||u_n||^2}
&=&\frac{1}{2}||w_n^+||^2-\frac{1}{2}||w^-_n||^2-\int
_{\mathbb{R}^{N}}\frac{F(x,u_n)}{||u_n||^2}dx \nonumber\\
&=&\frac{1}{2}||w_n^+||^2-\frac{1}{2}||w^-_n||^2-\int
_{\mathbb{R}^{N}}\frac{F_1(x,u_n)}{||u_n||^2}-\int
_{\mathbb{R}^{N}}\frac{F_2(x,u_n)}{||u_n||^2}dx,
\end{eqnarray}
where \begin{eqnarray}\label{m99xydhhd} F_1(x,t)=\left\{
\begin{array}
[c]{ll}
F(x,t),& |t|\leq\rho, \ x\in\mathbb{R}^N\\
0,& |t|>\rho, \ x\in\mathbb{R}^N
\end{array}
\right.\nonumber
\end{eqnarray}
and $F_2=F-F_1.$ From (\ref{nc99cidudyff}), we deduce that
$F_2\geq 0$ in $\mathbb{R}^N\times\mathbb{R}$. Let
\begin{eqnarray}\label{ji887u666999iu}
A:=\sup_{|t|\leq\rho, x\in\mathbb{R}^N}|f(x,t)/t|.
\end{eqnarray}
 From (\ref{l99ds7yydt}),
we have  $$A<\mu_0.$$  Then,
\begin{eqnarray}\label{bc88cydttd} |F_1(x,t)|\leq \frac{A}{2}t^2,\quad \forall (x,t)\in \mathbb{R}^N\times \mathbb{R}.
\end{eqnarray}    By (\ref{bc88cydttd}) and
(\ref{mmnc999ciu}), we have
\begin{eqnarray}\label{nc99cidudy}
\int
_{\mathbb{R}^{N}}\frac{|F_1(x,u_n)|}{||u_n||^2}dx\leq\frac{A}{2}||w_n||^2_{L^2}\leq\frac{A}{2\mu_0}||w_n||^2=
\frac{A}{2\mu_0}||w^+_n||^2+\frac{A}{2\mu_0}||w^-_n||^2.
\end{eqnarray}
Combining  (\ref{mmvn77dtdg}) and (\ref{nc99cidudy}), we get that
\begin{eqnarray}\label{mx99x0shgg}
-\frac{\mathcal{L}}{||u_n||^2}
\leq(\frac{1}{2}+\frac{A}{2\mu_0})||w_n^+||^2-(\frac{1}{2}-\frac{A}{2\mu_0})||w^-_n||^2-\int
_{\mathbb{R}^{N}}\frac{F_2(x,u_n)}{||u_n||^2}dx.
\end{eqnarray}

First, we consider the case $w\neq 0.$ From  $\lim_{|t|\rightarrow
\infty}f( x,t)/t = +\infty $ (see $\bf (f_2)$), we have
\begin{eqnarray}\label{nv99viuuvj}
\liminf_{|t|\rightarrow\infty}\frac{F_2(x,t)}{t^2}=+\infty.
\end{eqnarray}
 Note
that for $x\in\left\{  x\in\mathbb{R}^{N}\mid w\neq0\right\} $, we
have $| u_{n}(x)| \rightarrow+\infty$. Together with
(\ref{nv99viuuvj}), this implies that
\begin{eqnarray}\label{kgf888f9000}
\liminf_{n\rightarrow\infty}\int_{\{ x\in\mathbb{R}^{N}\ |\ w\neq0
\}}\frac{F_2(x,u_n)}{u_n^2}w^2_ndx=+\infty. \end{eqnarray} By
$F_2\geq 0,$ we get that
\begin{eqnarray}\label{n8888cyyy11} \int
_{\mathbb{R}^{N}}\frac{F_2(x,u_n)}{||u_n||^2}dx=\int
_{\mathbb{R}^{N}}\frac{F_2(x,u_n)}{u_n^2}w^2_ndx\geq \int _{\{
x\in\mathbb{R}^{N}\ |\  w\neq0 \}}\frac{F_2(x,u_n)}{u_n^2}w^2_ndx.
\end{eqnarray}  Combining  (\ref{kgf888f9000}) and (\ref{n8888cyyy11}) yields
\begin{eqnarray}\label{nnvb77dyd}
\liminf_{n\rightarrow\infty}\Big((\frac{1}{2}+\frac{A}{2\mu_0})||w_n^+||^2-(\frac{1}{2}-\frac{A}{2\mu_0})||w^-_n||^2
-\int
_{\mathbb{R}^{N}}\frac{F_2(x,u_n)}{||u_n||^2}dx\Big)=-\infty.\nonumber
\end{eqnarray}
It contradicts (\ref{mx99x0shgg}), since
$-\mathcal{L}/||u_n||^2\rightarrow0$ as $n\rightarrow\infty$.

Second, we consider the case $w= 0.$ In this case,
$\lim_{n\rightarrow\infty}||w^+_n||=0.$ It follows that
$$\lim_{n\rightarrow\infty}||w^-_n||=1,$$ since $||w_n||=1$ and
$||w_n||^2=||w^+_n||^2+||w^-_n||^2$. Therefore, the right hand
side of (\ref{mx99x0shgg}) is less than
$-\frac{1}{2}(\frac{1}{2}-\frac{A}{2\mu_0})$ when $n$ is large
enough. However, as $n\rightarrow\infty$, the left hand side of
(\ref{mx99x0shgg}) converges to zero. It also induces a
contradiction.

{\bf Step 3.} We shall prove that $\limsup_{ |||u|||\rightarrow
0}\Phi (u)\leq 0<\inf_N\Phi.$

From (\ref{bc88cydttd}) and the fact that $A<\mu_0$ and $F_2\geq
0$, we deduce that
\begin{eqnarray}\label{m00cudyyd}
\Phi(u) &=& \frac{1}{2}||u^+||^2-
\frac{1}{2}||u^-||^2-\int_{\mathbb{R}^N}F_1(x,u)dx-\int_{\mathbb{R}^N}F_2(x,u)dx\nonumber\\&\leq&
\frac{1}{2}||u^+||^2-
\frac{1}{2}||u^-||^2+\frac{A}{2}||u||^2_{L^2}-\int_{\mathbb{R}^N}F_2(x,u)dx\nonumber\\
&\leq&\frac{1}{2}||u^+||^2-
\frac{1}{2}||u^-||^2+\frac{A}{2\mu_0}||u||^2\nonumber\\
&=&(\frac{1}{2}+\frac{A}{2\mu_0})||u^+||^2-
(\frac{1}{2}-\frac{A}{2\mu_0})||u^-||^2\nonumber\\
&\leq&(\frac{1}{2}+\frac{A}{2\mu_0})||u^+||^2\leq(\frac{1}{2}+\frac{A}{2\mu_0})|||u|||^2.
\end{eqnarray}
This implies that $\limsup_{ |||u|||\rightarrow 0}\Phi (u)\leq
0<\inf_N \Phi.$

\medskip

Combining   Step 1-Step 3, we get the desired results of this
Lemma. \hfill$\Box$

\section{Boundedness of $(\overline{C})_c$ sequence  and proof of the main results}\label{nfh88fuy6fyf}

According to  Definition \ref{bcvujnbytg} in the appendix, a
sequence $\{u_n\}\subset X$ is called a $(\overline{C})_c$
sequence of $\Phi$ if
$$\sup_n\Phi(u_n)\leq c\quad \mbox{and}\quad
\lim_{n\rightarrow\infty}(1+||u_n||)||\Phi'(u_n)||_{X'}=0.$$
\begin{lemma}\label{1qx4esdftrg} Suppose that $\bf (v)$ and $\bf ( f_{1})-\bf ( f_{3}) $
  are satisfied. Let $\{u_n\}$ be a $(\overline{C})_c$ sequence of $\Phi$.
  Then
\begin{eqnarray}
\lim_{n\rightarrow\infty}\int_{\varpi_n}|u_n|^2dx=0,\
\lim_{n\rightarrow\infty}\int_{\varpi_n}|u_n|^pdx=0,\nonumber
\end{eqnarray}
where $\varpi_n=\{x\in\mathbb{R}^N\ |\ |u_n(x)|\geq M\}$ and $p$
and $M$ come from $\bf (f_1)$ and $\bf (f_3)$ respectively.
\end{lemma}
\noindent{\bf Proof.}   Let
$\widetilde{\varpi}^+_n=\{x\in\mathbb{R}^N\ |\ u_n(x)\geq M/2\}$
and $v_n=\max\{ u_n-M/2,0\}$. Then
\begin{eqnarray}\label{hhfnyydtrrss}
\int_{\mathbb{R}^N} |\nabla
v_n|^2dx=\int_{\widetilde{\varpi}^+_n}|\nabla
u_n|^2dx\leq\int_{\mathbb{R}^N} |\nabla u_n|^2dx,\
\int_{\mathbb{R}^N}  v_n^2dx=\int_{\widetilde{\varpi}^+_n}
v_n^2dx\leq\int_{\widetilde{\varpi}^+_n}
u_n^2dx\leq\int_{\mathbb{R}^N} u_n^2dx.\nonumber
\end{eqnarray}
It follows that $||v_n||=O( ||u_n||)$. Together with the fact that
$\{u_n\}$ is a $(\overline{C})_c$ sequence for $\Phi=-J,$ this
implies that
\begin{eqnarray}\label{jv99v8f7f7f}
o(1)=\langle\Phi'(u_n),v_n\rangle=-\langle J'(u_n),v_n\rangle.
\end{eqnarray}
By  (\ref{nc99ciduff}),
$$a:=\inf_{|t|\geq M/2,\ x\in\mathbb{R}^N}(\frac{f(x,t)}{t}-V_-(x))>0.$$
Then by (\ref{nv88ufyyfh1}), (\ref{jv99v8f7f7f}),
 and the fact that $u_n\geq v_n\geq 0$ on
$\widetilde{\varpi}^+_n$, we get that
\begin{eqnarray}\label{nv88fyftddd}
o(1)&=&\langle
J'(u_n),v_n\rangle\nonumber\\
&=&\int_{\mathbb{R}^N}(\nabla u_n\nabla
v_n+Vu_nv_n)dx+\int_{\mathbb{R}^N}f(x,u_n)v_ndx\nonumber\\
&=&\int_{\widetilde{\varpi}^+_n}|\nabla
v_n|^2dx+\int_{\widetilde{\varpi}^+_n}V_+(x)u_nv_ndx
+\int_{\widetilde{\varpi}^+_n}(\frac{f(x,u_n)}{u_n}-V_-)u_nv_ndx\nonumber\\
&\geq&\int_{\widetilde{\varpi}^+_n}|\nabla
v_n|^2dx+a\int_{\widetilde{\varpi}^+_n}v^2_ndx=\int_{\mathbb{R}^N}|\nabla
v_n|^2dx+a\int_{\mathbb{R}^N}v^2_ndx,
\end{eqnarray}
where $V_+=V+V_-\geq 0$ in $\mathbb{R}^N$.  Together with the
Sobolev inequality $\int_{\mathbb{R}^N}|\nabla
v_n|^2dx+a\int_{\mathbb{R}^N}v^2_ndx\geq
C(\int_{\mathbb{R}^N}|v_n|^pdx)^{2/p},$ this yields
\begin{eqnarray}\label{mcjudydte11}
\lim_{n\rightarrow\infty}\int_{\widetilde{\varpi}^+_n}|v_n|^2dx=\lim_{n\rightarrow\infty}\int_{\mathbb{R}^N}|v_n|^2dx=0,
\lim_{n\rightarrow\infty}\int_{\widetilde{\varpi}^+_n}|v_n|^pdx=\lim_{n\rightarrow\infty}\int_{\mathbb{R}^N}|v_n|^pdx=0.
\end{eqnarray}
Because
$$\varpi_n^+:=\{x\in\mathbb{R}^N\ |\ u_n(x)\geq M \}\subset\widetilde{\varpi}_n^+
$$
and  $v_n\geq u_n/2>0$ on $\varpi^+_n$, we get from
(\ref{mcjudydte11}) that
\begin{eqnarray}\label{bcyyfhtfgddd}
\lim_{n\rightarrow\infty}\int_{\varpi^+_n}|u_n|^2dx=0,\
\lim_{n\rightarrow\infty}\int_{\varpi^+_n}|u_n|^pdx=0.
\end{eqnarray}
Similarly, we can prove that
\begin{eqnarray}\label{tryrtrteteee}
\lim_{n\rightarrow\infty}\int_{\varpi^-_n}|u_n|^2dx=0,\
\lim_{n\rightarrow\infty}\int_{\varpi^-_n}|u_n|^pdx=0,
\end{eqnarray}
where $\varpi_n^-:=\{x\in\mathbb{R}^N\ |\ -u_n(x)\geq M \}.$ The
result of this lemma follows from (\ref{bcyyfhtfgddd}) and
(\ref{tryrtrteteee}).\hfill$\Box$

\begin{lemma}\label{x4esdftrg} Suppose that $\bf (v)$ and $\bf ( f_{1})-\bf ( f_{3}) $
  are satisfied. Let $\{u_n\}$ be a $(\overline{C})_c$ sequence of $\Phi$. Then
\begin{eqnarray}
\sup_{n}||u_n||<+\infty.\nonumber
\end{eqnarray}
\end{lemma}
\noindent{\bf Proof.} From
$(1+||u_n||)||\Phi'(u_n)||_{X'}\rightarrow 0$,   we get that
$\langle \Phi'(u_n), u^\pm_n\rangle=o(1)$. Then, by
(\ref{nc88c766f5}), we have
\begin{eqnarray}
||u^\pm_n||^2=\pm\int_{\mathbb{R}^N}f(x,u_n)u^\pm_ndx+o(1).\nonumber
\end{eqnarray}
It follows that
\begin{eqnarray}\label{hf55drzm}
||u_n||^2=\int_{\mathbb{R}^N}f(x,u_n)(u^+_n-u^-_n)dx+o(1).
\end{eqnarray}
From $\bf (f_1)$ and $\bf (f_2)$, we deduce that  there exists
$C_2>0$ such that
 \begin{eqnarray}\label{mn88fyfgged}
 |f(x,t)|\leq  |t|+C_2|t|^{p-1},\ \forall
 t\in\mathbb{R}.
 \end{eqnarray}
 Note that
$u^+$ and $u^-$ are orthogonal with respect to
$(\cdot,\cdot)_{L^2}.$ Then, by (\ref{mmnc999ciu}), we have
\begin{eqnarray}\label{cb88dyggdaa}
\int_{\mathbb{R}^N}|u^+-u^-|^2dx=\int_{\mathbb{R}^N}|u^+|^2dx+\int_{\mathbb{R}^N}|u^-|^2dx\leq\mu^{-1}_0||u||^2,\
\forall u\in X.\end{eqnarray} Let
\begin{eqnarray}\label{bcgttdydyyd}
D_1=\sup_{\rho\leq|t|\leq M,
x\in\mathbb{R}^N}\frac{f^2(x,t)}{|\widetilde{F}(x,t)|}
\end{eqnarray} and recall that $A=\sup_{|t|\leq\rho,
x\in\mathbb{R}^N}|f(x,t)/t|$ (see (\ref{ji887u666999iu})). Using
the H\"older inequality, from
 (\ref{hf55drzm}), (\ref{mn88fyfgged}) and (\ref{cb88dyggdaa}), we have
\begin{eqnarray}\label{jcnhhfbggdiij}
||u_n||^2&=&\Big(\int_{\{x\ |\ |u_n|\leq\rho\}}+\int_{\{x\ |\
\rho<|u_n|\leq M\}}+\int_{\{x\ |\ |u_n|>
M\}}\Big)f(x,u_n)(u^+_n-u^-_n)dx+o(1)\nonumber\\
&\leq&
A(\int_{\mathbb{R}^N}|u_n|^2dx)^{\frac{1}{2}}(\int_{\mathbb{R}^N}|u^+_n-u^-_n|^2dx)^{\frac{1}{2}}\nonumber\\
&&+(\int_{\{x\ |\ \rho<|u_n|\leq
M\}}f^2(x,u_n)dx)^{\frac{1}{2}}(\int_{\mathbb{R}^N}|u^+_n-u^-_n|^2dx)^{\frac{1}{2}}\nonumber\\
&&+(\int_{\{x\ |\ |u_n|>
M\}}|u_n|^2dx)^{\frac{1}{2}}(\int_{\mathbb{R}^N}|u^+_n-u^-_n|^2dx)^{\frac{1}{2}}\nonumber\\
&&+C_2(\int_{\{x\ |\ |u_n|>
M\}}|u_n|^pdx)^{\frac{p-1}{p}}(\int_{\mathbb{R}^N}|u^+_n-u^-_n|^pdx)^{\frac{1}{p}}\nonumber\\
&\leq&
A\mu^{-1}_0||u_n||^2+D^{\frac{1}{2}}_1\mu^{-\frac{1}{2}}_0(\int_{\{x\
|\ \rho<|u_n|\leq
M\}}|\widetilde{F}(x,u_n)|dx)^{\frac{1}{2}}||u_n||\nonumber\\
&&+\mu^{-\frac{1}{2}}_0(\int_{\{x\ |\ |u_n|>
M\}}|u_n|^2dx)^{\frac{1}{2}}||u_n||+C'C_2(\int_{\{x\ |\ |u_n|>
M\}}|u_n|^pdx)^{\frac{p-1}{p}}||u_n||,
\end{eqnarray}
where the positive constant $C'$ comes from the Sobolev inequality
$||u||_{L^p(\mathbb{R}^N)}\leq C||u||,$ $\forall u\in X.$  By
Lemma \ref{1qx4esdftrg}, we have
\begin{eqnarray}\label{cbgftyryyr88u}
\int_{\{x\ |\ |u_n(x)|\geq M\}}u^2_ndx=o(1),\ \int_{\{x\ |\
|u_n(x)|\geq M\}}|u_n|^pdx=o(1).
\end{eqnarray}
Combining (\ref{cbgftyryyr88u}) with (\ref{jcnhhfbggdiij}) yields
that
\begin{eqnarray}\label{nvoodjhhdg}
||u_n||^2\leq
A\mu^{-1}_0||u_n||^2+D^{\frac{1}{2}}_1\mu^{-\frac{1}{2}}_0(\int_{\{x\
|\ \rho<|u_n|\leq
M\}}|\widetilde{F}(x,u_n)|dx)^{\frac{1}{2}}||u_n||+o(||u_n||).
\end{eqnarray}

From $\sup_n\Phi(u_n)\leq c$ and
$(1+||u_n||)||\Phi'(u_n)||_{X'}\rightarrow 0$, we obtain
\begin{eqnarray}\label{bc66drdft2}
o(1)+c\geq\Phi(u_{n})-\frac{1}{2}\langle \Phi'(u_{n}),
u_n\rangle=\int_{\mathbb{R}^N}\widetilde{F}(x,u_n)dx.\nonumber
\end{eqnarray}
Together with (\ref{8r7bcvgftfudyff}), this implies
\begin{eqnarray}\label{fh99fijuufyf}
&&\int_{\{x\ |\ \rho<|u_n(x)|\leq
M\}}|\widetilde{F}(x,u_n)|dx\nonumber\\
&=&\int_{\{x\ |\
\rho<|u_n(x)|\leq M\}}\widetilde{F}(x,u_n)dx\nonumber\\
&\leq&- \int_{\{x\ |\ |u_n(x)|\leq \rho\}}\widetilde{F}(x,u_n)dx-
\int_{\{x\ |\ |u_n(x)|\geq
M\}}\widetilde{F}(x,u_n)dx+c+o(1)\nonumber\\
 &\leq& \int_{\{x\
|\ |u_n(x)|\leq \rho\}}|\widetilde{F}(x,u_n)|dx+ \int_{\{x\ |\
|u_n(x)|\geq M\}}|\widetilde{F}(x,u_n)|dx+|c|+o(1).
\end{eqnarray}
 From $\bf (f_1)$ and $\bf (f_2)$, we deduce that there exists
 $C_3>0$ such that
 \begin{eqnarray}\label{bfuu88fuyyf}
 |\widetilde{F}(x,t)|\leq |t|^2+C_3|t|^{p},\ \forall
 t\in\mathbb{R},
 \end{eqnarray}
Let $D_2=\sup_{|t|\leq\rho,
x\in\mathbb{R}^N}\frac{|\widetilde{F}(x,t)|}{t^2}$.  Combining
(\ref{bfuu88fuyyf}), (\ref{cbgftyryyr88u}) with
(\ref{fh99fijuufyf}) yields that
\begin{eqnarray}\label{nb8fdr5er}
&&\int_{\{x\ |\
\rho<|u_n(x)|\leq M\}}|\widetilde{F}(x,u_n)|dx\nonumber\\
&\leq& \int_{\{x\ |\ |u_n(x)|\leq \rho\}}|\widetilde{F}(x,u_n)|dx+
\int_{\{x\ |\ |u_n(x)|\geq
M\}}|\widetilde{F}(x,u_n)|dx+|c|+o(1)\nonumber\\
&\leq&D_2\int_{\{x\ |\ |u_n(x)|\leq \rho\}}u^2_ndx+\int_{\{x\ |\
|u_n(x)|\geq M\}}u^2_ndx+C_3\int_{\{x\ |\ |u_n(x)|\geq
M\}}|u_n|^pdx+|c|+o(1)\nonumber\\
&\leq&D_2\mu^{-1}_0||u_n||^2+\int_{\{x\ |\ |u_n(x)|\geq
M\}}u^2_ndx+C_3\int_{\{x\ |\
|u_n(x)|\geq M\}}|u_n|^pdx+|c|+o(1)\nonumber\\
&=&D_2\mu^{-1}_0||u_n||^2+|c|+o(1).\nonumber
\end{eqnarray}
Together with (\ref{nvoodjhhdg}), this implies
\begin{eqnarray}\label{nv99fdyfhhdg}
||u_n||^2\leq
\Big(A\mu^{-1}_0+D_1^{\frac{1}{2}}D_2^{\frac{1}{2}}\mu^{-1}_0\Big)
||u_n||^2+O(||u_n||).
\end{eqnarray}
From (\ref{l99ds7yydt}), we have
\begin{eqnarray}\label{nv999vivuuuv}
A\mu^{-1}_0+D_1^{\frac{1}{2}}D_2^{\frac{1}{2}}\mu^{-1}_0<1.
\end{eqnarray}
The boundedness of
 $\{||u_n||\}$  follows from (\ref{nv99fdyfhhdg}) and (\ref{nv999vivuuuv}) immediately.
\hfill$\Box$

\bigskip

\noindent{\bf Proof of Theorem \ref{th1}.} From the proof of Lemma
6.15 in \cite{willem}, we know that $\Phi'$ is weakly sequentially
continuous. Moreover, it is easy to see that $\sup_M
\Phi<+\infty.$ Then by Lemma \ref{bc77ftrg}, Lemma \ref{x4esdftrg}
and Theorem \ref{b99d443}, we deduce that there exists a bounded
$(\overline{C})_c$ sequence $\{u_n\}$ for $\Phi$ with
$c=\sup_M\Phi$ and $\inf_n|||u_n|||>0.$
   Up to a subsequence, either
\begin{itemize}
\item[{(i).}]
$\lim_{n\rightarrow\infty}\sup_{y\in\mathbb{R}^N}\int_{B_1(y)}|u_n|^2dx=0$,
or \item[{(ii).}] there exist $\varrho>0$ and $a_n\in\mathbb{Z}^N$
such that $\int_{B_1(a_n)}|u_n|^2dx\geq\varrho.$
\end{itemize}
If $(i)$ occurs, using the Lions lemma (see, for example,
\cite[Lemma 1.21]{willem}), a similar argument as for the proof of
\cite[Lemma 3.6]{szulkin} shows that
\begin{eqnarray}
\lim_{n\rightarrow\infty}\int_{\mathbb{R}^N}f(x,u_n)u^\pm_ndx=0.\nonumber
\end{eqnarray}
Then by (\ref{hf55drzm}), we have $||u_n||\rightarrow 0.$ This
contradicts $\inf_{n}|||u_n|||>0$. Therefore,  case $(i)$ cannot
occur. As case  $(ii)$ therefore occurs, $w_n=u_n(\cdot+a_n)$
satisfies $w_n\rightharpoonup u_0\neq 0$. From
$(1+||w_n||)||\Phi'(w_n)||_{X'}=(1+||u_n||)||\Phi'(u_n)||_{X'}\rightarrow
0$ and the weakly sequential continuity of $\Phi'$, we have that
$\Phi'(u_0)=0.$ Therefore, $u_0$ is a nontrivial solution of
Eq.(\ref{iiiuyhgt}). This completes the proof.\hfill$\Box$

\bigskip

\noindent{\bf Proof of Corollary \ref{nc88fyyyf2w}.} Let
$$f_\lambda(t)=|t|^{q-2}t-\lambda|t|^{p-2}t,\ t\in\mathbb{R}.$$ Since
$2<p<q<2^*,$  we deduce that $f_\lambda$ satisfies $\bf (f_1)$ and
$\bf (f_2).$ Because $q>p,$ there exists  $M>0$   such that
$$\inf_{|t|\geq M/2}(|t|^{q-2}-|t|^{p-2})>\max_{\mathbb{R}^N}V_-.$$
It follows that, for $0<\lambda\leq 1,$ $f_\lambda$ satisfies
(\ref{nc99ciduff}).

Let $$\rho=(q\lambda/p)^{1/(q-p)}.$$ Then, for $|t|\geq\rho,$
\begin{eqnarray}\label{nvciidyttd6}
F_\lambda(t)=\int^t_0f_\lambda(s)ds=\frac{1}{q}|t|^q-\frac{\lambda}{p}|t|^p\geq
0,
\end{eqnarray}
i.e., $F_\lambda$ satisfies (\ref{nc99cidudyff}).

 Let
$$\widetilde{F}_\lambda(t)=\frac{1}{2}tf_\lambda(t)-F_\lambda(t)
=(\frac{1}{2}-\frac{1}{q})|t|^q-\lambda(\frac{1}{2}-\frac{1}{p})|t|^p.$$
 If
$|t|\geq\rho,$ then
\begin{eqnarray}\label{nx99iuuuqa}
\widetilde{F}_\lambda(t)=(\frac{1}{2}-\frac{1}{q})|t|^q-\lambda(\frac{1}{2}-\frac{1}{p})|t|^p\geq\frac{q-p}{2q}|t|^q.
\end{eqnarray}
This shows that $\widetilde{F}_\lambda$ satisfies
(\ref{8r7bcvgftfudyff}).

It follows from \begin{eqnarray}\label{bv88dyfttdg}
0\leq|t|^{q-2}-\lambda|t|^{p-2}\leq |t|^{q-2} \ \mbox{ if}\
|t|\geq\rho \end{eqnarray} that
\begin{eqnarray}\label{nc99fuyyfh1}
f^2_\lambda(t)=(|t|^{q-2}-\lambda|t|^{p-2})^2t^2\leq
|t|^{2q-2}\quad\mbox{if}\quad |t|\geq\rho.
\end{eqnarray}
Let $\lambda$ be sufficiently small such that $M>\rho$. Combining
(\ref{nx99iuuuqa}) with (\ref{nc99fuyyfh1}) yields that
\begin{eqnarray}\label{bcvyysgtffda}
\frac{f^2_\lambda(t)}{\widetilde{F}_\lambda(t)}\leq\frac{2q}{q-p}|t|^{q-2}\leq\frac{2q}{q-p}M^{q-2}
\quad\mbox{if}\quad \rho\leq|t|\leq M.
\end{eqnarray}
Moreover, if $|t|\leq\rho$, we have
\begin{eqnarray}\label{bv99ttrfdee}
\Big|\frac{f_\lambda(t)}{t}\Big|=\Big||t|^{q-2}-\lambda|t|^{p-2}\Big|
\leq \rho^{q-2}+\lambda\rho^{p-2}
\end{eqnarray}
and
\begin{eqnarray}\label{nc88fyhgdffaaa}
\Big|\frac{\widetilde{F}_\lambda(t)}{t^2}\Big|\leq
(\frac{1}{2}-\frac{1}{q})|t|^{q-2}+\lambda(\frac{1}{2}-\frac{1}{p})|t|^{p-2}
\leq(\frac{1}{2}-\frac{1}{q})\rho^{q-2}+\lambda(\frac{1}{2}-\frac{1}{p})\rho^{p-2}.
\end{eqnarray}
Let $\lambda>0$ be sufficiently small such that
\begin{eqnarray}\label{nvb88yftftt}
\rho^{q-2}+\lambda\rho^{p-2}+\Big(\frac{2q}{q-p}M^{q-2}\Big)^{\frac{1}{2}}\cdot\Big((\frac{1}{2}-\frac{1}{q})\rho^{q-2}
+\lambda(\frac{1}{2}-\frac{1}{p})\rho^{p-2}\Big)^{\frac{1}{2}}<\mu_0.
\end{eqnarray}
It follows from $(\ref{bcvyysgtffda})-(\ref{nvb88yftftt})$ that
$f_\lambda$ and $\widetilde{F}_\lambda$ satisfy
(\ref{l99ds7yydt}). Therefore, we verified that $f_\lambda,$
$F_\lambda$ and $\widetilde{F}_\lambda$ satisfy
$(\ref{nc99ciduff})-(\ref{l99ds7yydt})$ if $\lambda>0$ is
sufficiently small. The result of this corollary follows from
Theorem \ref{th1} immediately.\hfill$\Box$

\section{Appendix}\label{nvb6ftrf}

In this section, we give a variant linking theorem which is a
generalization of the classical infinite-dimensional linking
theorem of \cite[Theorem 6.10]{willem} (see also \cite{KS}).

Before state this theorem, we give some notations and definitions.

Let $X$ be a separable Hilbert space with inner product
$(\cdot,\cdot)$ and norm $||\cdot||$, respectively. $X^\pm$ are
closed subspaces of $X$ and $X=X^+\oplus X^-.$ Let $\{e^-_k\}$ be
the total orthonormal sequence in $X^-$. Let
\begin{eqnarray}\label{nchu777d6d66}
Q:X\rightarrow X^+,\ P:X\rightarrow X^-
\end{eqnarray} be the
orthogonal projections. We define
\begin{eqnarray}\label{o88dytgdg}
|||u|||=\max\Big\{||Qu||,\sum^{\infty}_{k=1}\frac{1}{2^{k+1}}|(Pu,e^-_k)|\Big\}
\end{eqnarray}
on $X.$ Then $$||Qu||\leq|||u|||\leq ||u||,\ \forall u\in X.$$
Moreover, if $||u_n||$ is bounded and $|||u_n-u|||\rightarrow 0$,
then  $\{u_n\}$ weakly converges to $u$ in $X.$
 The topology
generated by $|||\cdot|||$ is denoted by $\tau$, and all
topological notations related to it will include this symbol.

\medskip

Let $R > r > 0$ and  $u_0 \in X^+$ with $||u_0||
 = 1$.
 Set
\begin{eqnarray}\label{ncg7777qwwq}
N = \{u \in X^+ \ |\  ||u||
 = r\},\  M = \{u+tu_0\ |\  u\in X^- ,\ t\geq 0,\
||u+tu_0|| \leq R\}. \end{eqnarray} Then, $M$ is a submanifold of
$X^- \oplus \mathbb{R}^+u_0$ with boundary
\begin{eqnarray}\label{m9155sfddd}
\partial M=\{u\in X^-\ |\ ||u||\leq
R\}\cup\{u+tu_0\ |\ u\in X^-,\ t>0,\ ||u+tu_0||=R \}.
\end{eqnarray}
\begin{definition}\label{bcvujnbytg}
Let  $\phi\in  C^1(X, \mathbb{R})$. A sequence $\{u_n\} \subset X$
is called a  $(\overline{C})_c$ sequence  for $\phi$, if $$
\sup_n\phi(u_n) \leq c\quad \mbox{ and}\quad
(1+||u_n
||)||\phi'(u_n)||_{X'}\rightarrow 0,\ \mbox{ as}\
n\rightarrow\infty.$$
\end{definition}

The main result of this section is  the following   theorem:
\begin{Theorem} \label{b99d443} If
 $H\in C^1(X,\mathbb{R})$ satisfies
\begin{description}  \item {$(a)$}  $H'$  is weakly
sequentially continuous, i.e., if $u\in X$ and $\{u_n\}\subset X$
are such that $u_n\rightharpoonup u$, then, for any $\varphi\in
X,$ $\langle H'(u_n),\varphi\rangle\rightarrow \langle
H'(u),\varphi\rangle$.
\item {$(b)$} There exist $\delta>0,$ $u_0
\in X^+  $ with $||u_0||=1$, and $R
> r
> 0$ such that
\begin{eqnarray}\label{nx99s8s7yy}
\inf_N H>\max\Big\{\sup_{\partial M} H,
\sup_{|||u|||\leq\delta}H(u)\Big\}
\end{eqnarray}
and
\begin{eqnarray}\label{m99iuxyyxy}
\sup_MH<+\infty,
\end{eqnarray}
\end{description} Then
there exists a  $(\overline{C})_c$  sequence $\{u_n\}$ for $H$
with $c= \sup_MH$ such that
\begin{eqnarray}\label{gft666f5fff11cv}
\inf_{n}|||u_n|||\geq\delta/2.
\end{eqnarray}
\end{Theorem}

\begin{remark}
In this theorem, the $\tau$-upper semi-continuous assumption of
$H$ in the classical infinite-dimensional linking theorem (see
\cite[Theorem 6.10]{willem}) is replaced by
$\inf_NH>\sup_{|||u|||\leq\delta}H(u)$ for some $\delta>0$ and the
result is replaced  with  existence of a $(\overline{C})_c$
sequence $\{u_n\}$ of $H$ satisfying (\ref{gft666f5fff11cv}).
\end{remark}

\bigskip

\noindent{\bf Proof of Theorem \ref{b99d443}.} Arguing indirectly,
assume that  the result does not hold. Then, there exists
$\epsilon>0$ such that
\begin{eqnarray}\label{nchyyftdg}
(1+||u||)||H'(u)||_{X'}\geq\epsilon,\ \forall u\in E
\end{eqnarray}
where $$E=\{u\in X\ |\ H(u)\leq d+\epsilon\}\cap\{u\in X\ |\
|||u|||\geq\delta/2\}$$ and $$d=\sup_{M}H.$$ From
(\ref{nx99s8s7yy}), we can choose $\epsilon$ such that
\begin{eqnarray}\label{nc88c7d6ttd}
0<\epsilon<\inf_NH-\max\Big\{\sup_{\partial M} H,
\sup_{|||u|||\leq\delta}H(u)\Big\}.
\end{eqnarray}

$\bf{Step 1.}$ A  vector field in a $\tau$-neighborhood of $E$.

Let
\begin{eqnarray}\label{hvcnbgdtdt66t11}
b=\inf_N H,\quad T=2(d-b+2\epsilon)/\epsilon,\quad R=(1+\sup_{u\in
M}||u||)e^T
\end{eqnarray}
and
\begin{eqnarray}\label{hhbc777cyttdg}
B_R=\{u\in X\ |\ ||u||\leq R\}.
\end{eqnarray}

For every $u\in E\cap B_R,$  there exists $\phi_u\in X$ with
$||\phi_u||=1$ such that $\langle
H'(u),\phi_u\rangle\geq\frac{3}{4}||H'(u)||_{X'}$. Then,
(\ref{nchyyftdg}) implies
\begin{eqnarray}\label{n99v8f77fy}
(1+||u||)\langle H'(u),\phi_u\rangle>\frac{1}{2}\epsilon.
\end{eqnarray}
From the definition of $|||\cdot|||$, we deduce that  if a
sequence $\{u_n\}\subset  E\cap B_R$  $\tau$-converges  to $u\in
X$, i.e., $|||u_n-u|||\rightarrow 0$, then $u_n\rightharpoonup u$
in $X$ (see Remark 6.1 of \cite{willem}).   By the weakly
sequential continuity of $H'$, we get that for any $\varphi\in X$,
$\langle H'(u_n),\varphi\rangle\rightarrow\langle
H'(u),\varphi\rangle$. This implies that $H'$ is
$\tau$-sequentially continuous in $E\cap B_R$. By
(\ref{n99v8f77fy}), the $\tau$-sequential continuity of $H'$ in
$E\cap B_R$ and the weakly lower semi-continuity of the norm
$||\cdot||$, we get that  there exists a $\tau$-open neighborhood
$V_u$ of $u$ such that
\begin{eqnarray}\label{nv88vuijf}
\langle H'(v),(1+||u||)\phi_u\rangle>\frac{1}{2}\epsilon,\ \forall
v\in V_u,
\end{eqnarray}
and
\begin{eqnarray}\label{nncbyyft55er}
||(1+||u||)\phi_u||=1+||u||\leq 2(1+||v||),\ \forall v\in V_u.
\end{eqnarray}

Because $B_R$ is a bounded convex closed set in the Hilbert space
$X$, $B_R$ is a $\tau$-closed set. Therefore, $X\setminus B_R$ is
a $\tau$-open set.

 The family
$$\mathcal{N}=\{V_u\ |\ u\in E\cap B_R\}\cup\{X\setminus B_R\}$$
is a $\tau$-open covering of $E$. Let
$$\mathcal{V}=\Big(\bigcup_{u\in E\cap B_R} V_u\Big)\bigcup (X\setminus B_R).$$
Then, $\mathcal{V}$ is a $\tau$-open neighborhood of $E.$

Since $\mathcal{V}$ is metric, hence paracompact, there exists a
local finite $\tau$-open covering $\mathcal{M}=\{M_i\ |\ i\in
\Lambda\}$ of $\mathcal{V}$ finer than $\mathcal{N}$. If
$M_i\subset V_{u_i}$ for some $u_i\in E$, we choose
$\varpi_i=(1+||u_i||)\phi_{u_i}$ and if $M_i\subset X\setminus
B_R$, we choose $\varpi_i=0$.
 Let
$\{\lambda_i\ |\ i\in I\}$ be a $\tau$-Lipschitz continuous
partition of unity subordinated to $\mathcal{M}$. And let
\begin{eqnarray}\label{ncb77dtdf}
\xi(u):=\sum_{i\in I}\lambda_i(u)\varpi_i,\ u\in \mathcal{V}.
\end{eqnarray}
Since the  $\tau$-open covering $\mathcal{M}$ of $\mathcal{V}$ is
local finite,  each $u\in\mathcal{V}$ belongs to only finite many
sets $M_i$. Therefore,  for every  $u\in\mathcal{V}$, the sum in
(\ref{ncb77dtdf}) is only a finite sum. It follows that, for any
$u\in \mathcal{V},$ there exist a $\tau$-open neighborhood
$U_u\subset\mathcal{V}$ of $u$ and $L_u>0$ such that
 $\xi(U_u)$ is
contained in a finite-dimensional subspace of $X$ and
\begin{eqnarray}\label{mm99ifuyy}
||\xi(v)-\xi(w)||\leq L_u|||v-w|||,\ \forall v,w\in U_u.
\end{eqnarray}
Moreover,  by the definition of $\xi$, (\ref{nv88vuijf}) and
(\ref{nncbyyft55er}), we get that, for every $u\in \mathcal{V},$
\begin{eqnarray}\label{oikjnhh7yyt}
||\xi(u)||\leq 1+||u||\quad \mbox{and}\quad \langle H'(u),
\xi(u)\rangle\geq 0
\end{eqnarray}
and for every $u\in E\cap B_R,$
\begin{eqnarray}\label{hhcbtgdtdr}
\langle H'(u), \xi(u)\rangle> \frac{1}{2}\epsilon.
\end{eqnarray}

$\bf{Step 2.} $  Let $\theta$ be a smooth  function satisfying
$0\leq \theta\leq 1$ in $\mathbb{R}$ and
$$\theta(t)=\left\{
\begin{array}
[c]{ll}
0 , & t\leq \frac{2\delta}{3},\\
1, & t\geq \delta.
\end{array}
\right.\label{e1}$$ Let
$$
\chi(u)=\left\{\begin{array} [c]{ll}
-\theta(|||u|||)\xi(u),&\ u\in\mathcal{V}\ ,  \\
0, & |||u||| \leq \frac{2\delta}{3}. \end{array}\right. \label{e1}
$$ Then, $\chi$ is a vector field defined in
$$\mathcal{W}=\mathcal{V}\cup\{u\in X\ |\ |||u|||<\delta\}.$$
It is a $\tau$-open neighborhood of $H^{d+\epsilon}\cup(X\setminus
B_R)$, where
$$H^{d+\epsilon}:=\{u\in X\ |\ H(u)\leq d+\epsilon\}.$$ From
(\ref{mm99ifuyy}), (\ref{hhcbtgdtdr}) and the definition of
$\chi$,  we deduce that the mapping $\chi$ satisfies that
\begin{description}
\item{$(\bf{a}).$} each $u\in \mathcal{W}$ has a $\tau$-open set
$V_u$ such that $\chi(V_u)$ is contained in a finite-dimensional
subspace of $X,$

\item{$(\bf{b}).$} for any $u\in \mathcal{W},$ there exist a
$\tau$-open neighborhood  $U_u$ of $u$ and $L'_u>0$ such that
\begin{eqnarray}\label{m000oiu76}
||\chi(v)-\chi(w)||\leq L'_u|||v-w|||,\ \forall v,w\in U_u.
\end{eqnarray}
This means that $\chi$ is locally Lipschitz continuous and
$\tau$-locally Lipschitz continuous,

\item{$(\bf{c}).$}
\begin{eqnarray}\label{bcvyyftfg} ||\chi(u)||\leq 1+||u||,\ \forall u\in
\mathcal{W},\end{eqnarray}

\item{$(\bf{d}).$}
\begin{eqnarray}\label{b7dyttdra} \langle
H'(u),\chi(u)\rangle\leq 0,\ \forall u\in \mathcal{W}.
\end{eqnarray} and
\begin{eqnarray}\label{b7dyttdratr} \langle
H'(u),\chi(u)\rangle<-\frac{1}{2}\epsilon,\ \forall u\in \{u\in E\
|\ |||u|||\geq\delta\}\cap B_R.
\end{eqnarray}
\end{description}

$\bf{Step 3.} $ From (\ref{m000oiu76}) and the fact that
$|||v|||\leq ||v||$, $\forall v\in X$, we have
$$
||\chi(v)-\chi(w)||\leq L'_u||v-w||,\ \forall v,w\in U_u. $$ This
implies that $\chi$ is a local Lipschitz mapping under the
$||\cdot||$ norm. Then by the standard theory of ordinary
differential equation in Banach space,  we deduce that the
following initial value problem
  \begin{eqnarray}\label{iijuytytyyy16}
\left\{
\begin{array}
[c]{ll}
\frac{d\eta}{dt}=\chi(\eta),\\
\eta(0,u)=u\in \mathcal{W}.
\end{array}
\right.
\end{eqnarray}
has a unique solution in $\mathcal{W}$ , denoted by $\eta(t, u)$,
with right maximal interval of existence $[0, T(u))$.
 Furthermore, using
(\ref{m000oiu76}) and the Gronwall inequality (see, for example,
Lemma 6.9 of \cite{willem}), the similar argument  as the proof of
$\bf c)$ in \cite[Lemma 6.8]{willem}  yields that
\begin{description}
\item{$(\bf{A}).$} $\eta$ is $\tau$-continuous, i.e., if
 $u_n\in \mathcal{W}$,  $u_0\in
\mathcal{W}$, $0\leq t_n<T(u_n)$ and $0\leq t_0<T(u_0)$ satisfy
$|||u_n-u_0|||\rightarrow 0$ and $t_n\rightarrow t_0$, then
$|||\eta(t_n,u_n)-\eta(t_0,u_0)|||\rightarrow 0.$
\end{description}

From (\ref{b7dyttdra}), we have $$\frac{d}{dt}H(\eta(t,u))=\langle
H'(\eta(t,u)),\eta(t,u)\rangle\leq0.$$  Therefore, $H$ is
non-increasing along the flow $\eta.$  It follows that
$\{\eta(t,u)\ |\ 0\leq t\leq T(u)\}\subset H^{d+\epsilon}$ if
$u\in H^{d+\epsilon},$ i.e., $H^{d+\epsilon}$ is an invariant set
of the flow $\eta.$ Then, (\ref{bcvyyftfg}) and Theorem 5.6.1 of
\cite{laksh}  implies  that, for any $ u\in H^{d+\epsilon},$
   $T(u)=+\infty$.

   \medskip

{\bf Step 4. } We shall prove that
\begin{eqnarray}\label{mxjvchhvc99987}
\{\eta(t, u)\ |\ 0\leq t\leq T,\ u\in
M\}\subset B_R.
\end{eqnarray}

Let $u\in H^{d+\epsilon}.$ By the result in Step 3,  we have
$T(u)=+\infty$ and
$$\eta(t,u)=u+\int^t_0 \chi(\eta(s,u))ds,\quad \forall t\in[0,+\infty).$$
Together with (\ref{bcvyyftfg}), this yields
\begin{eqnarray}\label{ufghyyft88u6}
||\eta(t,u)||\leq ||u||+\int^t_0 ||\chi(\eta(s,u))||ds\leq
||u||+\int^t_0(1+ ||\eta(s,u))||)ds.\nonumber
\end{eqnarray}
Then, by the Gronwall inequality (see, for example, Lemma 6.9 of
\cite{willem}), we get that
\begin{eqnarray}\label{jfuudt5555}
||\eta(t,u)||\leq (1+||u||)e^t-1, \quad \forall t\in[0,+\infty).
\end{eqnarray}
Since $M\subset H^{d+\epsilon}$, by (\ref{jfuudt5555}) and the
definition of $R$ (see \ref{hvcnbgdtdt66t11}), we get
(\ref{mxjvchhvc99987}).

\medskip

{\bf Step 5. } From  the choice of $\epsilon$ (see
(\ref{nc88c7d6ttd})), we have
$$\sup_{|||u|||\leq\delta}H<b-\epsilon.$$
It follows that
$$\{u\in X\ |\ |||u|||\leq\delta\}\subset H^{b-\epsilon}
:=\{u\in X\ |\ H(u)\leq b-\epsilon\}. $$ Together with
(\ref{b7dyttdratr}), this  yields
\begin{eqnarray}\label{hhdyttd}
\langle H'(u),\chi(u)\rangle<-\frac{1}{2}\epsilon,\ \forall u\in
H^{d+\epsilon}_{b-\epsilon}\cap B_R, \end{eqnarray} where
$$H^{d+\epsilon}_{b-\epsilon}:= \{u\in X\ |\ b-\epsilon\leq
H(u)\leq d+\epsilon\}.$$

We show that, for any $u\in M$, $H(\eta(T,u))\leq b-\epsilon.$
Arguing indirectly, assume that this were not true.  Then,  there
exists $u\in M$ such that $H(\eta(T,u))> b-\epsilon.$ Since $H$ is
non-increasing along the flow $\eta$, from (\ref{mxjvchhvc99987}),
we deduce that  $\{\eta(t,u)\ |\ 0\leq t\leq T\}\subset
H^{d+\epsilon}_{b-\epsilon}\cap B_R$. Then, by (\ref{hhdyttd}),
\begin{eqnarray}\label{nc88cudyd}
H(\eta(T,u))&=&H(\eta(0,u))+\int^T_0\Big\langle
H'(\eta(s,u)),\chi(\eta(s,u))\Big
\rangle ds\nonumber\\
&\leq&H(\eta(0,u))+\int^T_0(-\frac{1}{2}\epsilon)ds\nonumber\\
&\leq&d+\epsilon-\frac{1}{2}\epsilon T=b-\epsilon.
\end{eqnarray}
This contradicts $H(\eta(T,u))>b-\epsilon$. Therefore, we have
 \begin{description}
\item{$(\bf{B}).$} $\eta(T, M)\subset H^{b-\epsilon}.$
\end{description}
Moreover, using the result $\bf (a)$ in Step 2 and  and the fact
that $\eta$ is $\tau$-continuous (see $\bf (A)$), the similar
argument as the proof of the result $b)$  of \cite[Lemma
6.8]{willem} yields that
\begin{description}
\item{$(\bf{C}).$} Each point $(t,u)\in [0,T]\times
H^{d+\epsilon}$ has a $\tau$-neighborhood $N_{(t,u)}$ such that
$$\{v-\eta(s,v)\ |\ (s,v)\in N_{(t,u)}\cap ([0,T]\times
H^{d+\epsilon})\}$$ is contained in a finite-dimensional subspace
of $X.$
\end{description}

$\bf{Step 6.}$  Let $$h:[0,T]\times M\rightarrow X,\
h(t,u)=P\eta(t,u)+(||Q\eta(t,u)||-r)u_0$$ where $P, Q$, $r$ and
$u_0$ are defined in (\ref{nchu777d6d66}), (\ref{ncg7777qwwq}) and
(\ref{m9155sfddd}). Then $$0\in h(t,M)\Leftrightarrow
\eta(t,M)\cap N\neq\emptyset.$$ From $\inf_NH>\sup_{\partial M} H$
(see (\ref{nx99s8s7yy})) and the fact that, for any $u\in X$, the
function $H(\eta(\cdot, u))$ is non-increasing, we deduce that
$\inf_NH>\sup_{u\in\partial M}H(\eta(t,u))$, $\forall t\in [0,T]$.
Therefore,
\begin{eqnarray}\label{tgyyyuujpoik}0\not\in h(t,\partial M), \ \forall t\in
[0,T]. \end{eqnarray}  Since $\eta$ has the properties $\bf (A)$
and $\bf (C)$ obtained in step 3 and step 5 respectively and $h$
satisfies (\ref{tgyyyuujpoik}), there is an appropriate degree
theory for $\deg(h(t,\cdot), M,0)$ (see Proposition 6.4 and
Theorem 6.6 of \cite{willem}).  Then, the same argument as the
proof of Theorem 6.10 of \cite{willem} yields that
$$\deg(h(T,\cdot), M, 0)=\deg(h(0,\cdot), M, 0)\neq 0.$$
It follows that $0\in h(T,M)$ and $\eta(T,M)\cap N\neq\emptyset.$
Therefore, there exists $u\in M$ such that $H(\eta(T,u))\geq b.$
It contradicts the property $\bf (B)$ obtained  in step 5. This
completes the proof of this theorem.\hfill$\Box$

\end{document}